\documentclass[12pt]{article}

\usepackage{amsmath}
\usepackage{amssymb}
\usepackage{latexsym}
\usepackage{dsfont}
\begin{document}

\newtheorem{theo}{Theorem}
\newtheorem{prop}{Proposition}
\newtheorem{lemma}{Lemma}
\newtheorem{cor}{Corollary}
\newtheorem{defi}{Definition}
\newtheorem{exple}{Example}
\newtheorem{Rk}{Remark}

\begin{center}
{\Large \bf \sc Large and moderate deviations principles\\ for recursive kernel estimators of a
\vspace{0.25cm} 
\\ 
 multivariate density and its partial derivatives
}
\vspace{1.cm}

Abdelkader Mokkadem \hspace{1.2cm} Mariane Pelletier \hspace{1.2cm} Baba Thiam

\vspace{0.3cm}
(mokkadem, pelletier, thiam)@math.uvsq.fr

\vspace{0.5cm}

{\small Universit\'e de Versailles-Saint-Quentin\\
D\'epartement de Math\'ematiques\\
45, Avenue des Etats-Unis\\
78035 Versailles Cedex\\
France }
\vspace{1.cm}

{\bf Abstract:}\ \parbox[t]{9.cm}{\small
 In this paper we prove large and moderate deviations principles for the recursive kernel estimator of a probability density function and its partial derivatives. Unlike the density estimator, the derivatives estimators exhibit a quadratic behaviour not only for the moderate deviations scale but also for the large deviations one. We provide results both for the pointwise and the uniform deviations.}
\vspace{2.cm}

\rm

\end{center}

\noindent
{\bf Key words and phrases~:}\\
Multivariate recursive kernel estimation of a density and its derivatives~;
Large and moderate deviations principles
\vspace{0.5cm}

\newpage
\section{Introduction}
Let $X_1,\ldots,X_n$ be a sequence of independent and identically distributed $\mathbb{R}^d$-valued random vectors with bounded probability density $f$. Let $(h_n)$ be a positive sequence such that $\lim_{n \to \infty}h_n=0$ and $\sum_nh_n^d=\infty$; the recursive kernel estimator of $f$ is defined as
\begin{eqnarray}
\label{rec}
f_n(x)=\frac{1}{n}\sum_{i=1}^{n} \frac{1}{h_i^d}K\left(\frac{x-X_i}{h_i}\right)
\end{eqnarray}
where the kernel $K$ is a continuous function such that $\lim_{\|x\| \to +\infty}K(x)=0$ and $\int_{\mathbb{R}^d} K(x)dx=1$.
The estimate (\ref{rec}) is a recursive version of the well-known Rosenblatt kernel estimate (see Rosenblatt (1956) and Parzen (1962)); it was first discussed by Wolverton and Wagner (1969), Yamato (1971), and Davies (1973). The estimator (\ref{rec}) is easily updated each time an additionnal observation becomes available without resorting to past data, through the recursive relationship
\begin{eqnarray*}
f_n(x)=\frac{n-1}{n}f_{n-1}(x)+\frac{1}{nh_n^d}K\left(\frac{x-X_n}{h_n}\right).
\end{eqnarray*}  
 The weak and strong consistency of the recursive estimator of the density was studied by many authors; let us cite, among many others, Devroye (1979), Menon, Prasad and Singh (1984), and  Wertz (1985). The law of the iterated logarithm of the recursive density estimator was established by Wegman and Davies (1979) and Roussas (1992).  For other works on recursive density estimation, the reader is referred to the papers of Wegman (1972), Ahmad and Lin (1976), and Carroll (1976).\\

Recently, large and moderate deviations results have been proved for the Rosenblatt density estimator and its derivatives. The large deviations principle has been studied by Louani (1998) and Worms (2001). Gao (2003) and Mokkadem, Pelletier and Worms (2005) extend these results and provide moderate deviations principles. The large and moderate deviations of the derivatives of the Rosenblatt density estimator are given in Mokkadem, Pelletier and Worms (2005). The purpose of this paper is to establish large and moderate deviations principles for the recursive density estimator $f_n$ and its derivatives.\\

Let us recall that a $\mathbb{R}^m$-valued sequence $(Z_n)_{n \geq 1}$ satisfies a large deviations principle (LDP) with speed $(\nu_n)$ and  good rate function $I$ if:
\begin{description}
\item (a) $(\nu_n)$ is a positive sequence such that $(\nu_n)\uparrow \infty$;
\item (b) $I:\mathbb{R}^m \rightarrow [0,\infty]$ has compact level sets;
\item (c) for every borel set $B\subset \mathbb{R}^m$,
\end{description}
\begin{eqnarray*}
-\inf_{x \in \overset{\circ}{B}}I(x) & \leq & 
\liminf_{n \to \infty}\nu_n^{-1}\log\mathbb{P}\left[Z_n\in B\right]\\ & \leq &
\limsup_{n \to  \infty}\nu_n^{-1}\log\mathbb{P}\left[Z_n\in B\right]  \leq  
-\inf_{x \in \overline{B}}I(x),
\end{eqnarray*}
where $\overset{\circ}{B}$ and $\overline{B}$ denote the interior and the closure of $B$ respectively. Moreover, let $(v_n)$ be a nonrandom sequence that goes to infinity; if $(v_nZ_n)$ satisfies a LDP, then $(Z_n)$ is said to satisfy a moderate deviations principle (MDP).\\

For any $d$-uplet $\displaystyle[\alpha]=\Big(\alpha_1,\ldots,\alpha_d\Big)\in \mathbb{N}^d$, set $|\alpha|=\alpha_1+\cdots+\alpha_d$, let 
\begin{eqnarray*}
\partial^{[\alpha]}f(x)=\frac{\partial^{|\alpha|}f}{\partial x_1^{\alpha_1}\ldots\partial x_d^{\alpha_d}}(x)
\end{eqnarray*}
denote the $\displaystyle[\alpha]$-th partial derivative of $f$ (if $|\alpha|=0$, then $\partial^{[\alpha]}f \equiv f$) and, for any $j \in \mathbb{N}$, let $D^{(j)}f$ denote the $j$-th differential of $f$. The recursive kernel estimator of the $\displaystyle[\alpha]$-th partial derivative of $f$ is defined as 
\begin{eqnarray*}
\partial^{[\alpha]}f_n(x)=\frac{1}{n}\sum_{i=1}^{n} \frac{1}{h_i^{d+|\alpha|}}\partial^{[\alpha]}K\left(\frac{x-X_i}{h_i}\right),
\end{eqnarray*}
where the kernel $K$ is chosen such that $\partial^{[\alpha]}K \not\equiv 0$ and the bandwidth such that $\sum_nh_n^{d+2|\alpha|}=\infty$.\\

Our first aim is to establish pointwise LDP for the recursive kernel density estimator $f_n$. It turns out that expliciting the rate function in this case is more complex than either for the Rosenblatt kernel estimator, or for the derivatives estimators. That is the reason why, in this particular framework, we only consider bandwidths defined as $\left(h_n\right)\equiv \left(cn^{-a}\right)$ with $c>0$ and $a \in ]0,1/d[$. We then prove that the sequence  $\left(f_n(x)-f(x)\right)$ satisfies a LDP with speed $\left(\sum_{i=1}^nh_i^d\right)$ and rate function 
\begin{eqnarray*}
I_x:t\mapsto f(x)(1-ad)I\left(\frac{t}{f(x)(1-ad)}+\frac{1}{1-ad}\right)
\end{eqnarray*}
where $I(t)$ is the Fenchel-Legendre transform of the function $\psi$ defined as follows: 
\begin{eqnarray*}
\label{psi}
 \psi(u)  =  \int_{[0,1]\times\mathbb{R}^d}s^{-ad}\left(e^{s^{ad}\frac{u}{1-ad}K(z)}-1\right)dsdz.
\end{eqnarray*}
\\

Our second aim is to provide pointwise LDP for the derivative estimators $\partial^{[\alpha]}f_n$ (with $|\alpha| \geq 1$). In this case, we consider more general bandwidths defined as $h_n=h(n)$ for all $n$, where $h$ is a regularly varying function with exponent $(-a)$, $a \in ]0,1/\left(d+2|\alpha|\right)[$. We prove that the sequence
\begin{eqnarray*}
\left(\partial^{[\alpha]}f_n(x)-\partial^{[\alpha]}f(x)\right)
\end{eqnarray*}
satisfies a LDP of speed $\left(\sum_{i=1}^nh_i^{d+2|\alpha|}\right)$ and quadratic rate function  $J_{[\alpha]}:\mathbb{R}\rightarrow \mathbb{R}$ defined by
\begin{eqnarray}
\label{der}
\left\{
\begin{array}{ll}
\mbox{if} \quad f(x)  \neq  0, \quad J_{[\alpha],x}:t\mapsto \frac{t^2\left(1-a^2(d+2|\alpha|)^2\right)}{2f(x)\int_{\mathbb{R}^d}\left[\partial^{[\alpha]}K(z)\right]^2dz}\\
\mbox{if} \quad f(x)  =   0, \quad J_{[\alpha],x}(0)=0\quad \mbox{and}\quad J_{[\alpha],x}(t)=\infty \quad \mbox{for}\quad t\neq 0.
\end{array}\right.
\end{eqnarray}
\\

Our third aim is to prove pointwise MDP for the density estimator and for its derivatives. For any $d$-uplet $[\alpha]$ such that $|\alpha| \geq 0$, any positive sequence $(v_n)$ satisfying
\begin{eqnarray*}
\lim_{n \to \infty}v_n=\infty \quad \mbox{and}\quad \lim_{n \to \infty}\frac{v_n^2}{\sum_{i=1}^nh_i^{d+2|\alpha|}}=0,
\end{eqnarray*}
and general bandwidths $(h_n)$, we prove that the random sequence 
\begin{eqnarray*}
 v_n\left(\partial^{[\alpha]}f_n(x)-\partial^{[\alpha]}f(x)\right)
\end{eqnarray*}
satisfies a LDP of speed $\left(\sum_{i=1}^nh_i^{d+2|\alpha|}/v_n^2\right)$ and rate function $J_{[\alpha],x}(\cdot)$ defined by Equation \eqref{der}.
\\

Finally, we give a uniform version of the previous results. More precisely, let $U$ be a subset of $\mathbb{R}^d$; we establish large and moderate deviations principles for the sequence \\$\left(\sup_{x \in U}\left|\partial^{[\alpha]}f_n(x)-\partial^{[\alpha]}f(x)\right|\right)$ in the case either $U$ is bounded or all the moments of $f$ are finite.
\section{Assumptions and Results}
\subsection{Pointwise LDP for the density estimator}
The assumptions required on the kernel $K$ and the bandwidth $(h_n)$ are the following.
\begin{description}
\item (H1) $K:\mathbb{R}^d\rightarrow \mathbb{R}$ is a bounded and integrable function , $\int_{\mathbb{R}^d}K(z)dz=1$ and $\lim_{\|z\| \to \infty}K(z)=0$.
\item (H2) $h_n=cn^{-a}$ with  $0< a < 1/d$ and $c >0$.               
\end{description}
Before stating our results, we need to introduce the rate function for the LDP of the density estimator.
Let $\psi:\mathbb{R}\rightarrow \mathbb{R}$ and $I:\mathbb{R}\rightarrow \mathbb{R}$ be the functions defined as:
\begin{eqnarray*}
 \psi(u)  =  \int_{[0,1]\times\mathbb{R}^d}s^{-ad}\left(e^{s^{ad}\frac{u}{1-ad}K(z)}-1\right)dsdz\quad \mbox{and}\quad
I(t)  =  \sup_{u \in \mathbb{R}}\left\{ut-\psi(u)\right\}
\end{eqnarray*}
(where $s \in [0,1]$, $z \in \mathbb{R}^d$) and set 
\begin{eqnarray*}
S_+ =\left\{x \in \mathbb{R}^d ; K(x)>0\right\} \quad \mbox{and}\quad S_- = \left\{x \in \mathbb{R}^d ; K(x)<0\right\}. 
\end{eqnarray*}
The following proposition gives the properties of the functions $\psi$ and $I$; in particular, the behaviour of the rate function $I$, which differs depending on whether $K$ is non-negative or not, is explicited.
\begin{prop}\label{rate}
Let $\lambda$ be the Lebesgue measure on $\mathbb{R}^d$ and let Assumption (H1) hold.
\begin{description}
\item (i) $\psi$ is strictly convex, twice continuously differentiable on $\mathbb{R}$, and $I$ is a good rate function on $\mathbb{R}$.
\item (ii) If $\lambda(S_-)=0$, then $I(t)=+\infty$ when $t<0$, $\displaystyle I(0)=\lambda(S_+)/(1-ad)$, $I$ is strictly convex on $\mathbb{R}$ and continuous on $]0,+\infty[$, and for any $t>0$
\begin{eqnarray}
\label{mel}
I(t)=t\left(\psi'\right)^{-1}(t)-\psi\left(\left(\psi'\right)^{-1}(t)\right).
\end{eqnarray}
\item (iii) If $\lambda (S_-)>0$, then $I$ is finite and strictly convex on $\mathbb{R}$ and \eqref{mel} holds for any $t \in \mathbb{R}$.
\item (iv) In both cases, the strict minimum  of $I$ is achieved by  $I(1/(1-ad))=0$.  
\end{description}
\end{prop}

\paragraph{Remark}
The following relations are straightforward, and will be used in the sequel:
\begin{eqnarray}
\label{cre}
I(t) & = & 
\left\{
\begin{array}{ll}
\sup_{u>0}\{ut-\psi(u)\} \quad \mbox{if}\quad t>1/(1-ad)\\
\sup_{u<0}\{ut-\psi(u)\} \quad \mbox{if}\quad t<1/(1-ad).
\end{array}\right.
\end{eqnarray}
\\
We can now state the LDP for the density estimator.
\begin{theo}\label{fat}\textbf{(Pointwise LDP for the density estimator)}$ $ \\
Let Hypotheses (H1)-(H2) hold  and assume that $f$ is continuous at $x$. Then, the sequence $\left(f_n(x)-f(x)\right)$ satisfies a LDP with speed $(\sum_{i=1}^nh_i^d)$ and rate function defined as follows:
\begin{eqnarray*}
\left\{
\begin{array}{ll}
\mbox{if} \quad f(x)  \neq  0, \quad I_x:t\mapsto f(x)(1-ad)I\left(\frac{1}{1-ad}+\frac{t}{f(x)(1-ad)}\right)\\
\mbox{if} \quad f(x)  =   0, \quad I_x(0)=0\quad \mbox{and}\quad I_x(t)=+\infty \quad \mbox{for}\quad t\neq 0.
\end{array}\right.
\end{eqnarray*}

\end{theo}

\subsection{Pointwise LDP for the derivatives estimators}
Let $[\alpha]$ be a $d$-uplet such that $|\alpha|\geq 1$. To establish pointwise LDP for $\partial^{[\alpha]}f_n$, we need the following assumptions.
\begin{description}
\item (H3) $h_n=h(n)$ where the function $h$ is locally bounded and varies regularly with exponent $(-a)$, $0<a<1/(d+2|\alpha|)$.
\item (H4) i) $K$ is $|\alpha|$-times differentiable on $\mathbb{R}^d$ and $\lim_{\|x\| \to \infty}\|D^{(j)}K(x)\|=0$ for any $j\in \left\{0,\ldots,|\alpha|-1\right\}$. \\
         ii) $\partial^{[\alpha]}K :\mathbb{R}^d\rightarrow \mathbb{R}$ is a bounded and integrable function and $\int_{\mathbb{R}^d}\left[\partial^{[\alpha]}K(x)\right]^2dx \neq 0$.
\item (H5) $f$ is $|\alpha|$-times differentiable on $\mathbb{R}^d$ and its $j$-th differentials $D^{(j)}f$ are bounded on $\mathbb{R}^d$ for any $j\in \left \{0,\ldots,|\alpha|-1\right\}$.       
\end{description}

\paragraph{Remark} 
A positive (not necessarily monotone) function $L$ defined on $]0,\infty[$ is slowly varying if $\lim_{t \to \infty}L(tx)/L(t)=1$; a function $G$ is said to vary regularly with exponent $\rho$, $\rho \in \mathbb{R}$, if and only if there exists a slowly varying function $L$ such that $G(x)=x^{\rho}L(x)$ (see, for example, Feller (1970) page 275). Typical examples of regularly varying functions with exponent $\rho$ are $x^\rho$, $x^\rho\log x$, $x^\rho\log \log x$, $\displaystyle x^\rho\log x/\log \log x$, and so on. An important consequence of (H3) which will be used in the sequel is:
\begin{eqnarray}
\label{rel}
\mbox{if} \quad\beta a <1,\quad \mbox{then} \quad\lim_{n \to \infty}\frac{1}{nh_n^\beta}\sum_{i=1}^nh_i^{\beta}= \frac{1}{1-a\beta}.
\end{eqnarray}
\begin{theo}\label{ldp}\textbf{(Pointwise LDP for the derivatives estimators)}$ $ \\
Let $|\alpha|\geq 1$ and assume that (H1), (H3)-(H5) hold and that $\partial^{[\alpha]}f$ is continuous at $x$. Then, the sequence $\left(\partial^{[\alpha]}f_n(x)-\partial^{[\alpha]}f(x)\right)$ satisfies a LDP with speed $(\sum_{i=1}^nh_i^{d+2|\alpha|})$ and rate function $J_{[\alpha],x}$ defined by \eqref{der}.
\end{theo}
\subsection{Pointwise MDP for the density estimator and its derivatives}
Let $(v_n)$ be a positive sequence; we assume that
\begin{description}
\item (H6) $\lim_{n \to \infty}v_n=\infty$ and $\displaystyle\lim_{n \to \infty}\frac{v_n^2}{\sum_{i=1}^nh_i^{d+2|\alpha|}}=0$.\\
          
\item (H7) i) There exists an integer $q \geq 2$ such that $\forall s \in \{1,\ldots,q-1\},\forall j \in \{1,\ldots,d\}$, $\displaystyle\int_{\mathbb{R}^d}y_j^sK(y)dy_j=0$, and $\displaystyle\int_{\mathbb{R}^d}\big|y_j^qK(y)\big|dy<\infty$.\\
          ii) $\displaystyle\lim_{n \to \infty}\frac{v_n}{n}\sum_{i=1}^nh_i^q=0$.\\        iii) $\partial^{[\alpha]}f$ is $q$-times differentiable on $\mathbb{R}^d$ and $M_q = \sup_{x \in \mathbb{R}^d}\|D^q\partial^{[\alpha]}f(x)\|<+\infty$.
\end{description}
\paragraph{Remark} 
When $h_n=O(n^{-a})$, with $0<a<1/\left(d+2|\alpha|\right)$, (H6) and (H7)ii) hold for instance for $(v_n)\equiv (n^b)$ for any $b \in ]0,\min\{aq ; \left(1-a(d+2|\alpha|)\right)/2\}[$.
\\
The following theorem gives the MDP for the density estimator and its derivatives.
\begin{theo}\label{mdp}\textbf{(Pointwise MDP)} $ $ \\
When $|\alpha|=0$, let Assumptions (H1), (H2), (H6) and (H7) hold; when $|\alpha| \geq 1$, let (H1), (H3)-(H7) hold. If $\partial^{[\alpha]}f$ is q-times differentiable at $x$, then the sequence $\left(v_n\left(\partial^{[\alpha]}f_n(x)-\partial^{[\alpha]}f(x)\right)\right)$ satisfies a LDP with speed $\left(\sum_{i=1}^nh_i^{d+2|\alpha|}/v_n^2\right)$ and good rate function $J_{[\alpha],x}$ defined in (\ref{der}).
\end{theo}

\subsection{Uniform LDP and MDP for the density estimator and its derivatives}
 To establish uniform large deviations principles for the density estimator and its derivatives, we need the following additionnal assumptions:
\begin{description}
\item (H8) i) There exists $\xi>0$ such that $\int_{\mathbb{R}^d}\|x\|^{\xi}f(x)dx<\infty$.\\
          ii) $f$ is uniformly continuous.
\item (H9) i) $\partial^{[\alpha]}K$ is H\"older continuous.\\
          ii) There exists $\gamma>0$ such that $z \mapsto\|z\|^{\gamma}\partial^{[\alpha]}K(z)$ is a bounded function.
\item (H10) $\displaystyle\lim_{n \to \infty}\frac{v_n^2\log (1/h_n)}{\sum_{i=1}^nh_i^{d+2|\alpha|}}=0$ and $\displaystyle\lim_{n \to \infty}\frac{v_n^2\log v_n}{\sum_{i=1}^nh_i^{d+2|\alpha|}}=0$.
\item (H11) i) There exists $\zeta>0$ such that $\int_{\mathbb{R}^d}\|z\|^{\zeta}\left|K(z)\right|dz<\infty$.\\
           ii) There exists $\eta>0$ such that $z \mapsto\|z\|^{\eta}\partial^{[\alpha]}f(z)$ is a bounded function. 
\end{description}
\paragraph{Remark} 
When $h_n=O\left(n^{-a}\right)$ with $a \in ]0,1/(d+2|\alpha|)[$, (H10) holds for instance with $(v_n)\equiv (n^{b})$ for any $b \in ]0,\left(1-a(d+2|\alpha|)\right)/2[$.
\\

Set $U\subseteq \mathbb{R}^d$; in order to state  in a compact form the uniform large and moderate deviations principles for the density estimator and its derivatives on $U$, we consider the large deviations case as the special case $(v_n)\equiv 1$ and we set: 
\begin{eqnarray*}
g_{U}(\delta) & = &
\left\{
\begin{array}{l}
\|f\|_{U,\infty}(1-ad)I\left(\frac{1}{1-ad}+\frac{\delta}{\|f\|_{U,\infty}(1-ad)}\right) \quad \mbox{if}\quad |\alpha|=0\quad\mbox{and}\quad (v_n) \equiv 1\\
 \frac{\delta^2\left(1-a^2(d+2|\alpha|)^2\right)}{2\|f\|_{U,\infty}\int_{\mathbb{R}^d}\left[\partial^{[\alpha]}K\right]^2(z)dz}\quad \mbox{otherwise},
\end{array}\right.\\
\tilde{g}_U(\delta) & = & \min\{g_U(\delta),g_U(-\delta)\},
\end{eqnarray*}
where
$\|f\|_{U,\infty}=\sup_{x \in U}|f(x)|$.\\
\paragraph{Remark}
The functions $g_U(\cdot)$ and $\tilde{g}_U(\cdot)$ are non-negative, continuous, increasing on $]0,+\infty[$ and decreasing on $]-\infty,0[$, with a unique global minimum in $0$ $(\tilde{g}_U(0)=g_U(0)=0)$. They are thus good rate functions (and $g_U(\cdot)$ is strictly convex).\\
 
Theorem \ref{unidevldp} below states uniform LDP and MDP for $\left(\partial^{[\alpha]}f_n-\partial^{[\alpha]}f\right)$ on $U$ in the case $U$ is bounded, and Theorem \ref{uldp} in the case $U$ is unbounded.
\begin{theo}\label{unidevldp} \textbf{(Uniform deviations on a bounded set)}$ $ \\
In the case $|\alpha|=0$, let (H1), (H2), (H7), (H9)i), and (H10) hold. In the case $|\alpha|\geq 1$, let (H3)-(H5), (H7), (H9)i) and (H10) hold. Moreover, assume either that $(v_n)\equiv 1$ or that $(v_n)$ satisfies (H6). Then 
for any bounded subset $U$ of $\mathbb{R}^d$ and for all $\delta >0$,
\begin{eqnarray}
\label{bound}
\lim_{n \to \infty}\frac{v_n^2}{\sum_{i=1}^nh_i^{d+2|\alpha|}}\log\mathbb{P}\left[\sup_{x \in U}v_n\left|\partial^{[\alpha]}f_n(x)-\partial^{[\alpha]}f(x)\right|\geq \delta\right] & = & -\tilde{g}_U(\delta).
\end{eqnarray}
\end{theo}

\begin{theo}\label{uldp}\textbf{(Uniform deviations on an unbounded set)} $ $ \\
Let Assumptions (H1), (H7)-(H11) hold. Moreover, 
\begin{itemize}
\item in the case $|\alpha|=0$ and $(v_n)\equiv 1$, let (H2) hold;
\item in the case $|\alpha|\geq 1$ and $(v_n)\equiv 1$, or $|\alpha|\geq 0$ and $(v_n)$ satisfies (H6), let (H3)-(H5) hold.
\end{itemize}
Then for any subset $U$ of $\mathbb{R}^d$ and for all $\delta >0$,
\begin{eqnarray*}
-\tilde{g}_U(\delta) & \leq & \liminf_{n \to \infty}\frac{v_n^2}{\sum_{i=1}^nh_i^{d+2|\alpha|}}\log\mathbb{P}\left[\sup_{x \in U}v_n\left|\partial^{[\alpha]}f_n(x)-\partial^{[\alpha]}f(x)\right|\geq \delta\right]\\ 
& \leq & 
\limsup_{n \to \infty}\frac{v_n^2}{\sum_{i=1}^nh_i^{d+2|\alpha|}}\log\mathbb{P}\left[\sup_{x \in U}v_n\left|\partial^{[\alpha]}f_n(x)-\partial^{[\alpha]}f(x)\right|\geq \delta\right] \leq 
-\frac{\xi}{\xi+d}\tilde{g}_U(\delta).
\end{eqnarray*}
\end{theo}
The following corollary is a straightforward consequence of Theorem \ref{uldp}.
\begin{cor}\label{cor} $ $ \\
Under the assumptions of Theorem \ref{uldp}, if $\int_{\mathbb{R}^d}\|x\|^{\xi}f(x)dx<\infty$ $\forall \xi \in \mathbb{R}$, then for any subset $U$ of $\mathbb{R}^d$,
\begin{eqnarray}
\label{fo}
\lim_{n \to \infty}\frac{v_n^2}{\sum_{i=1}^nh_i^{d+2|\alpha|}}\log\mathbb{P}\left[\sup_{x \in U}v_n\left|\partial^{[\alpha]}f_n(x)-\partial^{[\alpha]}f(x)\right|\geq \delta\right] & = & -\tilde{g}_U(\delta).
\end{eqnarray}
\end{cor}
\paragraph{Comment}
Theorem \ref{unidevldp} and Corollary \ref{cor} are LDP for the sequence $\left(\sup_{x \in U}\left|f_n(x)-f(x)\right|\right)$. As a matter of fact, since the sequence $\left(\sup_{x \in U}\left|f_n(x)-f(x)\right|\right)$ is positive and since $\tilde{g}_U$ is continuous on $[0,+\infty[$, increasing and goes to infinity as $\delta \to \infty$, the application of Lemma $5$ in Worms (2001) allows to deduce from \eqref{bound} or \eqref{fo} that $\left(\sup_{x \in U}\left|f_n(x)-f(x)\right|\right)$ satisfies a LDP with speed $\left(\sum_{i=1}^nh_i^d\right)$ and good rate function $\tilde{g}_U$ on $\mathbb{R}_+$.
\section{Proofs}
Let $(\Psi_n^{[\alpha]})$ and $(B_n^{[\alpha]})$ be the sequences defined as
\begin{eqnarray*}
\Psi_n^{[\alpha]}(x) & = & \frac{1}{n}\sum_{i=1}^n\frac{1}{h_i^{d+|\alpha|}}\left(\partial^{[\alpha]}K\left(\frac{x-X_i}{h_i}\right)-\mathbb{E}\left[\partial^{[\alpha]}K\left(\frac{x-X_i}{h_i}\right)\right]\right),\\ B_n^{[\alpha]}(x) & = & \mathbb{E}\left[\partial^{[\alpha]}f_n(x)\right]-\partial^{[\alpha]}f(x). 
\end{eqnarray*}
We have:
\begin{eqnarray}
\label{dec}
\partial^{[\alpha]}f_n(x)-\partial^{[\alpha]}f(x) = \Psi_n^{[\alpha]}(x)+B_n^{[\alpha]}(x).
\end{eqnarray}
Theorems \ref{fat}, \ref{ldp} and \ref{mdp} are consequences of \eqref{dec} and the following propositions. 
\begin{prop}\label{vardp}\textbf{(Pointwise LDP and MDP for \boldmath$(\Psi_n^{[\alpha]})$)}
\begin{enumerate}
\item Under the assumptions (H1) and (H2), the sequence $\left(f_n(x)-\mathbb{E}\left(f_n(x)\right)\right)$ satisfies a LDP with speed $\left(\sum_{i=1}^nh_i^d\right)$ and rate function $I_x$.
\item Let $|\alpha|\geq 1$ and assume that (H3), (H4) hold, then the sequence $\left(\Psi_n^{[\alpha]}(x)\right)$ satisfies a LDP with speed $\left(\sum_{i=1}^nh_i^{d+2|\alpha|}\right)$ and rate function $J_{[\alpha],x}$.
\item When $|\alpha|=0$, let Assumptions (H1), (H2) and (H6) hold and when $|\alpha| \geq 1$, let (H3), (H4) and (H6) hold, then the sequence $\left(v_n\Psi_n^{[\alpha]}(x)\right)$ satisfies a LDP with speed $\left(\sum_{i=1}^nh_i^{d+2|\alpha|}/v_n^2\right)$ and rate function $J_{[\alpha],x}$. 
\end{enumerate}
\end{prop}

\begin{prop}\label{univardev}\textbf{(Uniform LDP and MDP for \boldmath$(\Psi_n^{[\alpha]})$)} 
\begin{enumerate}
\item In the case $|\alpha|=0$, let (H1), (H2), (H9)i) and (H10)  hold. In the case $|\alpha|\geq 1$, let (H3)-(H5), (H9)i) and (H10) hold. Moreover, assume either that $(v_n)\equiv1$ or that $(v_n)$ satisfies (H6); then 
for any bounded subset $U$ of $\mathbb{R}^d$ and for all $\delta >0$,
\begin{eqnarray*}
\lim_{n \to \infty}\frac{v_n^2}{\sum_{i=1}^nh_i^{d+2|\alpha|}}\log\mathbb{P}\left[\sup_{x \in U}v_n\left|\partial^{[\alpha]}\Psi_n(x)\right|\geq \delta\right] & = & -\tilde{g}_U(\delta).
\end{eqnarray*}

\item Let Assumptions (H1), (H8)-(H11) hold . Moreover,
\begin{itemize}
\item in the case $|\alpha|=0$ and $(v_n)\equiv 1$, let (H2) holds,
\item in the case either $|\alpha| \geq 1$ and $(v_n)\equiv 1$, or $|\alpha| \geq 0$ and $(v_n)$ satisfies (H6), let (H3)-(H5) hold.
\end{itemize}
Then for any subset $U$ of $\mathbb{R}^d$ and for all $\delta >0$,
\begin{eqnarray*}
-\tilde{g}_U(\delta) & \leq & \liminf_{n \to \infty}\frac{v_n^2}{\sum_{i=1}^nh_i^{d+2|\alpha|}}\log\mathbb{P}\left[\sup_{x \in U}v_n\left|\partial^{[\alpha]}\Psi_n(x)\right|\geq \delta\right]\\ 
& \leq & 
\limsup_{n \to \infty}\frac{v_n^2}{\sum_{i=1}^nh_i^{d+2|\alpha|}}\log\mathbb{P}\left[\sup_{x \in U}v_n\left|\partial^{[\alpha]}\Psi_n(x)\right|\geq \delta\right] \leq 
-\frac{\xi}{\xi+d}\tilde{g}_U(\delta).
\end{eqnarray*}
\end{enumerate}
\end{prop}

\begin{prop}\label{bia} \textbf{(Pointwise and uniform convergence rate of \boldmath $B_n^{[\alpha]}$)}\\
Let Assumptions (H1), (H3)-(H5) and (H7)i) hold. 
\begin{itemize}
\item [1)] If $\partial^{[\alpha]}f$ is $q$-times differentiable at $x$, then
\begin{eqnarray*}
\mathbb{E}\left(\partial^{[\alpha]}f_n(x)\right)-\partial^{[\alpha]}f(x)  & = & 
O\left(\frac{\sum_{i=1}^nh_i^q}{n}\right).
\end{eqnarray*}
\item [2)]If (H7)iii) holds, then:
\begin{eqnarray*}
\lim_{n \to \infty}\frac{n}{\sum_{i=1}^nh_i^q}\sup_{x \in \mathbb{R}^d}\left|\mathbb{E}\left(\partial^{[\alpha]}f_n(x)\right)-\partial^{[\alpha]}f(x)\right|  & \leq & \frac{M_q}{q!}\int_{\mathbb{R}^d}\|z\|^q\left|K(z)\right|dz.
\end{eqnarray*}
\end{itemize}
\end{prop}
Set $x \in \mathbb{R}^d$; since the assumptions of Theorems \ref{fat} and \ref{ldp} guarantee that $\lim_{n \to \infty}B_n^{[\alpha]}(x)=0$, Theorem \ref{fat} (respectively Theorem \ref{ldp}) is a straightforward consequence of the application of Part 1 (respectively of Part 2) of Proposition \ref{vardp}.
Moreover, under the assumptions of Theorem \ref{mdp}, we have, by application of Part 1 of Proposition \ref{bia}, $\lim_{n \to \infty}v_nB_n^{[\alpha]}(x)=0$; Theorem \ref{mdp} thus straightfully follows from the application of Part 3 of Proposition \ref{vardp}. Finally, Theorems \ref{unidevldp} and \ref{uldp} are obtained by applying Proposition \ref{univardev} together with the second part of Proposition \ref{bia}.\\ 

We now state a preliminary lemma, which will be used in the proof of Propositions \ref{vardp} and \ref{univardev}. For any $u \in \mathbb{R}$, set   
\begin{eqnarray*}
\Lambda_{n,x}(u) & = & 
\frac{v_n^2}{\sum_{i=1}^nh_i^{d+2|\alpha|}}\log\mathbb{E}\left[\exp\left(u\frac{\sum_{i=1}^nh_i^{d+2|\alpha|}}{v_n}\Psi_n^{[\alpha]}(x)\right)\right],\\
\Lambda_x^L(u)
& = & f(x)(1-ad)\left(\psi(u)-\frac{u}{1-ad}\right),\\
\Lambda_x^M(u) & = & \frac{u^2}{2(1-a^2\left(d+2|\alpha|)^2\right)}f(x)\int_{\mathbb{R}^d}\left[\partial^{[\alpha]}K(z)\right]^2dz. 
\end{eqnarray*}
\begin{lemma}\label{con}\textbf{(Convergence of \boldmath $\Lambda_{n,x}$)} $ $
\begin{itemize}
\item In the case $|\alpha|=0$ and $(v_n)\equiv 1$, let (H1) and (H2) hold;
\item In the case either $|\alpha|\geq 1$ and $(v_n)\equiv 1$, or $|\alpha|\geq 0$ and $(v_n)$ satisfies (H6), let (H1), (H3) and (H4) hold.
\end{itemize}
\begin{enumerate}
\item (Pointwise convergence)\\
If $f$ is continuous at $x$, then for all $u \in \mathbb{R}$,
\begin{eqnarray}
\label{lam}
\lim_{n \to \infty}\Lambda_{n,x}(u)  =  \Lambda_x(u) 
\end{eqnarray}

where 
\begin{eqnarray*}
 \Lambda_x(u) =  \left\{
\begin{array}{ll}
  \Lambda_x^{L}(u)\quad \mbox{when}\quad v_n \equiv 1\quad \mbox{and}\quad |\alpha|=0\\  \Lambda_x^{M}(u)\quad \mbox{when}\quad either \quad v_n \to \infty \quad \mbox{and}\quad |\alpha|\geq 0 \quad  \mbox{or} \quad v_n\equiv 1\quad \mbox{and}\quad|\alpha|\geq 1.  
\end{array}
\right.
\end{eqnarray*}

\item (Uniform convergence)\\
If $f$ is uniformly continuous, then the convergence (\ref{lam}) holds uniformly in $x\in U$.
\end{enumerate}
\end{lemma}
Our proofs are now organized as follows: Lemma \ref{con} is proved in Section 3.1, Proposition \ref{vardp} in Section 3.2, Proposition \ref{univardev} in Section 3.3 and Proposition \ref{bia} in Section 3.4. Section 3.5 is devoted to the proof of Proposition \ref{rate} on the rate function $I$.
\subsection{Proof of Lemma \ref{con}}
Set $u \in \mathbb{R}$, $Y_i=\partial^{[\alpha]}K\left(\frac{x-X_i}{h_i}\right)$ and $a_n=\sum_{i=1}^nh_i^{d+2|\alpha|}$. We have:
\begin{eqnarray*}
\Lambda_{n,x}(u) & = & \frac{v_n^2}{a_n}\log\mathbb{E}\left[\exp\left(u\frac{a_n}{v_n}\Psi_n^{[\alpha]}(x)\right)\right]\\ & = &
\frac{v_n^2}{a_n}\log\mathbb{E}\left[\exp\left(u\frac{a_n}{nv_n}\sum_{i=1}^n\frac{1}{h_i^{d+|\alpha|}}\left(Y_i-\mathbb{E}(Y_i)\right)\right)\right] \\ & = & 
\frac{v_n^2}{a_n}\sum_{i=1}^n\log\mathbb{E}\left[\exp\left(\frac{ua_n}{nv_nh_i^{d+|\alpha|}}Y_i\right)\right]-\frac{uv_n}{n}\sum_{i=1}^n\frac{1}{h_i^{d+|\alpha|}}\mathbb{E}(Y_i).
\end{eqnarray*}
By Taylor expansion, there exists $c_{i,n}$ between 1 and $\mathbb{E}\left[\exp\left(\frac{ua_n}{nv_nh_i^{d+|\alpha|}}Y_i\right)\right]$ such that 
\begin{eqnarray*}
\lefteqn{\log\left(\mathbb{E}\left[\exp\left(\frac{ua_n}{nv_nh_i^{d+|\alpha|}}Y_i\right)\right]\right)} \\ \nonumber  & = & 
\mathbb{E}\left[\exp\left(\frac{ua_n}{nv_nh_i^{d+|\alpha|}}Y_i\right)-1\right] -\frac{1}{2c_{i,n}^2}\left(\mathbb{E}\left[\exp\left(\frac{ua_n}{nv_nh_i^{d+|\alpha|}}Y_i\right)-1\right]\right)^2
\end{eqnarray*}
and $\Lambda_{n,x}$ can be rewritten as
\begin{eqnarray}
\label{titi}
\lefteqn{\Lambda_{n,x}(u)}\nonumber\\ & = & 
\frac{v_n^2}{a_n}\sum_{i=1}^n\mathbb{E}\left[\exp\left(\frac{ua_n}{nv_nh_i^{d+|\alpha|}}Y_i\right)-1\right]-\frac{v_n^2}{2a_n}\sum_{i=1}^n\frac{1}{c_{i,n}^2}\left(\mathbb{E}\left[\exp\left(\frac{ua_n}{nv_nh_i^{d+|\alpha|}}Y_i\right)\right]-1\right)^2 \nonumber\\ \mbox{} &&-\frac{uv_n}{n}\sum_{i=1}^n\frac{1}{h_i^{d+|\alpha|}}\mathbb{E}(Y_i).
\end{eqnarray}
For proving Lemma \ref{con}, we consider two cases:
\subsubsection{ First case: either \boldmath $v_n\to \infty$ or \boldmath $|\alpha|\geq 1$.} 
A Taylor's expansion implies the existence of  $c'_{i,n}$ between $0$ and $\displaystyle\frac{ua_n}{nv_nh_i^{d+|\alpha|}}Y_i$ such that 
\begin{eqnarray*}
\lefteqn{\mathbb{E}\left[\exp\left(\frac{ua_n}{nv_nh_i^{d+|\alpha|}}Y_i\right)-1\right]}\\ & = &
 \frac{ua_n}{nv_nh_i^{d+|\alpha|}}\mathbb{E}(Y_i)+\frac{1}{2}\left(\frac{ua_n}{nv_nh_i^{d+|\alpha|}}\right)^2\mathbb{E}(Y_i^2) + \frac{1}{6}\left(\frac{ua_n}{nv_nh_i^{d+|\alpha|}}\right)^3\mathbb{E}\left(Y_i^3e^{c'_{i,n}}\right).
\end{eqnarray*}
Therefore, 
\begin{eqnarray*}
\Lambda_{n,x}(u) & = & 
\frac{u^2a_n}{2n^2}\sum_{i=1}^n\frac{1}{h_i^{2d+2|\alpha|}}\mathbb{E}\left[Y_i^2\right]+R_{n,x}^{(1)}(u)\\ & = & 
 f(x)\frac{u^2a_n}{2n^2}\sum_{i=1}^n\frac{1}{h_i^{d+2|\alpha|}}\int_{\mathbb{R}^d}\left(\partial^{[\alpha]}K(z)\right)^2dz+R_{n,x}^{(1)}(u)+R_{n,x}^{(2)}(u)
\end{eqnarray*}
with 
\begin{eqnarray*}
R_{n,x}^{(1)}(u) & = & \frac{1}{6}\frac{u^3a_n^2}{n^3v_n}\sum_{i=1}^nh_i^{-3d-3|\alpha|}\mathbb{E}\left(Y_i^3e^{c'_{i,n}}\right)   -\frac{v_n^2}{2a_n}\sum_{i=1}^n\frac{1}{c_{i,n}^2}\left(\mathbb{E}\left[\exp\left(\frac{ua_n}{nv_nh_i^{d+|\alpha|}}Y_i\right)-1\right]\right)^2\\
 R_{n,x}^{(2)}(u) & = & \frac{u^2a_n}{2n^2}\sum_{i=1}^n\frac{1}{h_i^{d+2|\alpha|}}\int_{\mathbb{R}^d}\left(\partial^{[\alpha]}K(z)\right)^2\left[f(x-h_iz)-f(x)\right]dz.
\end{eqnarray*}
Using \eqref{rel}, one can show that
\begin{eqnarray*} 
\left|\frac{ua_n}{nv_nh_i^{d+|\alpha|}}Y_i\right|
& \leq & c_1\frac{|u|}{1-a\left(d+2|\alpha|\right)}\|\partial^{[\alpha]}K\|_{\infty}
\end{eqnarray*}
where $c_1$ is a positive constant and thus
\begin{eqnarray}
\label{fot}
c'_{i,n} & \leq & c_1 \frac{|u|}{1-a(d+2|\alpha|)}\|\partial^{[\alpha]}K\|_{\infty}
\end{eqnarray}
and
\begin{eqnarray} 
\label{fet}
\frac{1}{c_{i,n}^2} \leq \exp\left(2c_1\frac{|u|}{1-a(d+2|\alpha|)}\|\partial^{[\alpha]}K\|_{\infty}\right).
\end{eqnarray}
Noting that $\mathbb{E}|Y_i|^3 \leq h_i^d\|f\|_{\infty}\int_{\mathbb{R}^d}\left|\partial^{[\alpha]}K(z)\right|^3dz$ and using \eqref{rel} again and \eqref{fot}, there exists a positive constant $c_2$ such that, for $n$ large enough,
\begin{eqnarray}
\label{mami}
\lefteqn{\left|\frac{u^3a_n^2}{n^3v_n}\sum_{i=1}^nh_i^{-3d-3|\alpha|}\mathbb{E}\left(Y_i^3e^{c'_{i,n}}\right)\right|} \nonumber 
 \\ & \leq & 
c_2\frac{|u|^3e^{\frac{c_1|u|\|\partial^{[\alpha]}K\|_{\infty}}{1-a(d+2|\alpha|)}}h_n^{|\alpha|}}{\left(1+a(2d+3|\alpha|)\right)\left(1-a(d+2|\alpha|)\right)^2v_n}\|f\|_{\infty}\int_{\mathbb{R}^d}\left|\partial^{[\alpha]}K(z)\right|^3dz
\end{eqnarray}
which goes to $0$ as $n \to \infty$ since either $v_n \to \infty$ or $|\alpha| \geq 1$.\\
In the same way, in view of \eqref{rel} and \eqref{fet}, there exists a positive constant $c_3$ such that, for $n$ large enough,  
\begin{eqnarray}
\label{fit}
\lefteqn{\left|\frac{v_n^2}{2a_n}\sum_{i=1}^n\frac{1}{c_{i,n}^2}\left(\mathbb{E}\left[\exp\left(\frac{ua_n}{nv_nh_i^{d+|\alpha|}}Y_i\right)-1\right]\right)^2\right|}\nonumber \\ 
& \leq &c_3\|f\|_{\infty}^2e^{\frac{4c_1|u|}{1-a(d+2|\alpha|)}\|\partial^{[\alpha]}K\|_{\infty}}\left(\int_{\mathbb{R}^d}\left|\partial^{[\alpha]}K(z)\right|dz\right)^2h_n^d
\end{eqnarray}
which goes to $0$ as $n \to \infty$. The combination of \eqref{mami} and \eqref{fit} ensures that  $$\lim_{n \to \infty}\sup_{x \in \mathbb{R}^d}|R_{n,x}^{(1)}(u)|=0.$$
Now, since $f$ is continuous, we have $\lim_{i \to \infty}\left|f(x-h_iz)-f(x)\right|=0$, and thus, by the dominated convergence theorem, (H4)ii) implies that $\lim_{i \to \infty}\int_{\mathbb{R}^d}\left(\partial^{[\alpha]}K(z)\right)^2\left|f(x-h_iz)-f(x)\right|dz=0$.  
Since, in view of \eqref{rel}, the sequence $\displaystyle\left(\frac{u^2a_n}{2n^2}\sum_{i=1}^n\frac{1}{h_i^{d+2|\alpha|}}\right)$ is bounded, it follows that $\displaystyle\lim_{n \to \infty}|R_{n,x}^{(2)}(u)|=0$. The pointwise convergence \eqref{lam} then follows.\\
In the case $f$ is uniformly continuous, we have $\lim_{i \to \infty}\sup_{x \in \mathbb{R}^d}\left|f(x-h_iz)-f(x)\right|=0$ and thus, using the same arguments as previously, we obtain $\lim_{n \to \infty}\sup_{x \in \mathbb{R}^d}|R_{n,x}^{(2)}(u)|=0$.\\
We then deduce that $\lim_{n \to \infty}\sup_{x \in U}\left|\Lambda_{n,x}(u)-\Lambda_x(u)\right| =  0$ which concludes the proof of Lemma \ref{con} in this case.
\subsubsection{Second case: \boldmath$|\alpha|=0$ and \boldmath$(v_n)\equiv 1$.} 
 Using assumption (H2) and in view of (\ref{titi}), there exists $c>0$ such that
\begin{eqnarray*}
\Lambda_{n,x}(u) & = & \frac{1}{a_n}\sum_{i=1}^n\mathbb{E}\left[e^{\frac{ua_n}{nh_i^d}Y_i}-1\right]-\frac{1}{2a_n}\sum_{i=1}^n\frac{1}{c_{i,n}^2}\left(\mathbb{E}\left[e^{\frac{ua_n}{nh_i^d}Y_i}-1\right]\right)^2-\frac{u}{n}\sum_{i=1}^nh_i^{-d}\mathbb{E}(Y_i)\\ & = & 
\frac{1}{a_n}\sum_{i=1}^nh_i^d\int_{\mathbb{R}^d}\left[e^{\frac{cu}{1-ad}\left(\frac{i}{n}\right)^{ad}K(z)}-1-uK(z)\right]f(x)dz-R_{n,x}^{(3)}(u)+R_{n,x}^{(4)}(u)
\end{eqnarray*}
with
\begin{eqnarray*}
R_{n,x}^{(3)}(u) & = & \frac{1}{2a_n}\sum_{i=1}^n\frac{1}{c_{i,n}^2}\left(\mathbb{E}\left[e^{\frac{ua_n}{nh_i^d}Y_i}-1\right]\right)^2 \\
R_{n,x}^{(4)}(u)
 & = & \frac{1}{a_n}\sum_{i=1}^nh_i^d\int_{\mathbb{R}^d}\left(e^{\frac{ua_n}{nh_i^d}K(z)}-1\right)\left[f(x-h_iz)-f(x)\right]dz \\ \mbox{} & & -\frac{u}{n}\sum_{i=1}^n\int_{\mathbb{R}^d}K(z)\left[f(x-h_iz)-f(x)\right]dz.
\end{eqnarray*}
Using the bound \eqref{fit}, we have $\lim_{n \to \infty}\sup_{x \in \mathbb{R}^d}|R_{n,x}^{(3)}(u)|=0$.\\
 Since $|e^t-1|\leq |t|e^{|t|}$, we have
\begin{eqnarray*}
\left|R_{n,x}^{(4)}(u)\right|
 & \leq & \frac{1}{a_n}\sum_{i=1}^nh_i^d\int_{\mathbb{R}^d}\left|\left(e^{\frac{ua_n}{nh_i^d}K(z)}-1\right)\left[f(x-h_iz)-f(x)\right]\right|dz \\ \mbox{} & & +\frac{|u|}{n}\sum_{i=1}^n\int_{\mathbb{R}^d}\left|K(z)\right|\left|f(x-h_iz)-f(x)\right| \\ & \leq & 
\frac{|u|}{n}e^{\frac{c|u|}{1-ad}\|K\|_{\infty}}\sum_{i=1}^n\int_{\mathbb{R}^d}\left|K(z)\right|\left|f(x-h_iz)-f(x)\right|dz \\ \mbox{} & & + \frac{|u|}{n}\sum_{i=1}^n\int_{\mathbb{R}^d}\left|K(z)\right|\left|f(x-h_iz)-f(x)\right| \\ & = & \left(e^{\frac{c|u|}{1-ad}\|K\|_{\infty}}+1\right) \frac{|u|}{n}\sum_{i=1}^n\int_{\mathbb{R}^d}\left|K(z)\right|\left|f(x-h_iz)-f(x)\right|dz.
\end{eqnarray*}
In the case $f$ is continuous, we have $\lim_{n \to \infty}R_{n,x}^{(4)}(u)=0$ by the dominated convergence theorem.\\
In the case $f$ is uniformly continuous, set $\varepsilon >0$ and let $M>0$ such that $2\|f\|_{\infty}\int_{\|z\|\geq M}\left|K(z)\right|dz<\varepsilon/2$. We need to prove that for $n$ sufficiently large 
\begin{eqnarray*}
\sup_{x \in \mathbb{R}^d}\int_{\|z\| \leq M}\left|K(z)\right|\left|f(x-h_iz)-f(x)\right|dz & \leq & \varepsilon/2
\end{eqnarray*}
which is a straightforward consequence of the uniform continuity of $f$. 
It follows from analysis considerations that 
\begin{eqnarray*}
\lim_{n \to \infty}\Lambda_{n,x}(u) & = & f(x)\int_{\mathbb{R}^d}(1-ad)\left[\int_{0}^1s^{-ad}\left(e^{s^{ad}\frac{u}{1-ad}K(z)}-1-uK(z)\right)ds\right]dz\\ & = & 
\frac{f(x)}{\int_0^1s^{-ad}ds}\int_{0}^1s^{-ad}\left[\int_{\mathbb{R}^d}\left(e^{s^{ad}\frac{u}{1-ad}K(z)}-1-uK(z)\right)dz\right]ds.
\end{eqnarray*}
and thus Lemma \ref{con} is proved. $\square$
\subsection{Proof of Proposition \ref{vardp}}
To prove Proposition \ref{vardp}, we apply Proposition \ref{rate}, Lemma \ref{con} and the following result (see Puhalskii (1994)). 
\begin{lemma}\label{puh}
Let $(Z_n)$ be a sequence of real random variables, $(\nu_n)$ a positive sequence satisfying $\lim_{n \to \infty}\nu_n=+\infty$, and suppose that there exists some convex non-negative function $\Gamma$ defined (\emph{i.e.} finite) on $\mathbb{R}$ such that 
\begin{eqnarray*}
\forall u\in \mathbb{R}, \quad \lim_{n \to \infty}\frac{1}{\nu_n}\log\mathbb{E}\left[\exp(u\nu_nZ_n)\right] = \Gamma(u).
\end{eqnarray*}
If the Legendre transform $\tilde{\Gamma}$ of $\Gamma$ is a strictly convex function, then the sequence $(Z_n)$ satisfies a LDP of speed $(\nu_n)$ and good rate function $\tilde{\Gamma}$. 
\end{lemma}
In our framework, when $|\alpha|=0$ and $v_n\equiv 1$, we take $Z_n=f_n(x)-\mathbb{E}(f_n(x))$, $\nu_n=\sum_{i=1}^nh_i^d$ with $h_n=cn^{-a}$ where $0<a<1/d$ and $\Gamma=\Lambda_x^{L}$. In this case, the Legendre transform of $\Gamma=\Lambda_x^L$ is the rate function $I_x:t\mapsto f(x)(1-ad)I\left(\frac{t}{f(x)(1-ad)}+\frac{1}{1-ad}\right)$ which is strictly convex by Proposition \ref{rate}. Otherwise, we take $Z_n=v_n\left(\partial^{[\alpha]}f_n(x)-\mathbb{E}\left(\partial^{[\alpha]}f_n(x)\right)\right)$, $\nu_n=\sum_{i=1}^nh_i^{d+2|\alpha|}/v_n^2$ and $\Gamma=\Lambda_x^{M}$; $\tilde{\Gamma}$ is then the quadratic rate function $J_{[\alpha],x}$ defined in (\ref{der}) and thus  Proposition \ref{vardp} follows. $\square$
\subsection{Proof of Proposition \ref{univardev}}
In order to prove Proposition \ref{univardev}, we first establish some lemmas.

\begin{lemma}\label{sup}
Let $\phi:\mathbb{R}+\to \mathbb{R}$ be the function defined for $\delta>0$ as
\begin{eqnarray*}
\phi(\delta) & = & 
\left\{
\begin{array}{l}
(\psi')^{-1}\left(\frac{\delta}{\|f\|_{U,\infty}(1-ad)}+\frac{1}{1-ad}\right) \quad \mbox{if}\quad (v_n) \equiv 1\quad \mbox{and}\quad |\alpha|=0,\\
 \frac{\delta\left(1-a^2(d+2|\alpha|)^2\right)}{\|f\|_{U,\infty}\int_{\mathbb{R}^d}\left[\partial^{[\alpha]}K\right]^2(z)dz}\quad \mbox{otherwise}.
\end{array}\right.
\end{eqnarray*}
\begin{enumerate}
\item $\sup_{u\in \mathbb{R}}\{u\delta-\sup_{x \in U}\Lambda_x(u)\}$ equals $g_U(\delta)$ and is achieved for $u=\phi(\delta)>0$.
\item $\sup_{u\in \mathbb{R}}\{-u\delta-\sup_{x \in U}\Lambda_x(u)\}$ equals $g_U(-\delta)$ and is achieved for $u=\phi(-\delta)<0$.
\end{enumerate}
\end{lemma}
\subsubsection*{Proof of Lemma \ref{sup}}
We just prove the first part, the proof of the second one being similar. First, let us consider the case $(v_n)\equiv 1$ and $|\alpha|=0$. Since $e^t\geq 1+t$ $(\forall t)$, we have $\psi(u) \geq u/(1-ad)$ and therefore,
\begin{eqnarray*}
u\delta-\sup_{x \in U}\Lambda_x(u) & = & u\delta-\|f\|_{U,\infty}(1-ad)\left(\psi(u)-\frac{u}{1-ad}\right) \\ & = & 
\|f\|_{U,\infty}(1-ad)\left[u\left(\frac{\delta}{\|f\|_{U,\infty}(1-ad)}+\frac{1}{1-ad}\right)-\psi(u)\right]
\end{eqnarray*}
The function $u\mapsto u\delta-\sup_{x \in U}\Lambda_x(u)$ has second derivative $-\|f\|_{U,\infty}(1-ad)\psi''(u)<0$ and thus it has a unique maximum achieved for 
\begin{eqnarray*}
u_0=(\psi')^{-1}\left(\frac{\delta}{\|f\|_{U,\infty}(1-ad)}+\frac{1}{1-ad}\right).
\end{eqnarray*}
Now, since $\psi'$ is increasing and since $\displaystyle\psi'(0)=1/(1-ad)$, we deduce that $u_0>0$.\\
In the case $\lim_{n \to \infty}v_n=\infty$, Lemma \ref{sup} is established in the same way by noting that 
\begin{eqnarray*}
u\delta-\sup_{x \in U}\Lambda_x(u) = u\delta-\sup_{x \in U}\Lambda_x^M(u) = u\delta -\frac{u^2}{2\left(1-a^2\left(d+2|\alpha|\right)^2\right)}\|f\|_{U,\infty}\int_{\mathbb{R}^d}\left[\partial^{[\alpha]}K(z)\right]^2dz. \quad\square
\end{eqnarray*}

\begin{lemma}\label{sort} $ $
\begin{itemize}
\item In the case $|\alpha|=0$ and $(v_n)\equiv1$, let (H1) and (H2) hold;
\item In the case either $|\alpha|\geq 1$ and $(v_n)\equiv 1$, or $|\alpha|\geq 0$ and $(v_n)$ satisfies (H6), let (H1), (H3) and (H4) hold.
\end{itemize}
Then for any $\delta>0$,
\begin{eqnarray*}
\lim_{n \to \infty}\frac{v_n^2}{\sum_{i=1}^nh_i^{d+2|\alpha|}}\log\sup_{x \in U}\mathbb{P}\left[v_n\left(\partial^{[\alpha]}f_n(x)-\mathbb{E}\left(\partial^{[\alpha]}f_n(x)\right)\right)\geq \delta \right] & = & -g_U(\delta)\\
\lim_{n \to \infty}\frac{v_n^2}{\sum_{i=1}^nh_i^{d+2|\alpha|}}\log\sup_{x \in U}\mathbb{P}\left[v_n\left(\partial^{[\alpha]}f_n(x)-\mathbb{E}\left(\partial^{[\alpha]}f_n(x)\right)\right)\leq -\delta \right] & = & -g_U(-\delta)\\
\lim_{n \to \infty}\frac{v_n^2}{\sum_{i=1}^nh_i^{d+2|\alpha|}}\log\sup_{x \in U}\mathbb{P}\left[v_n\left|\partial^{[\alpha]}f_n(x)-\mathbb{E}\left(\partial^{[\alpha]}f_n(x)\right)\right|\geq \delta \right] & = & -\tilde{g}_U(\delta).
\end{eqnarray*}

\end{lemma}
\subsubsection*{Proof of Lemma \ref{sort}}
Set $\displaystyle b_n=\sum_{i=1}^nh_i^{d+2|\alpha|}/v_n^2$, $S_n(x)=v_n\Psi_n^{[\alpha]}(x)$, and $\delta>0$. In the sequel, $\Lambda_x(u)$ is defined as in (\ref{lam}).\\
We first note that, for any $u>0$,
\begin{eqnarray*}
\mathbb{P}\left[S_n(x)\geq \delta\right] & = & \mathbb{P}\left[e^{b_nuS_n(x)}\geq e^{b_nu\delta}\right]\\ & \leq & 
e^{-b_nu\delta}\mathbb{E}\left[e^{b_nuS_n(x)}\right]\\
& \leq & 
e^{-b_nu\delta}e^{b_n\Lambda_{n,x}(u)}\\
& \leq & 
e^{-b_n\left(u\delta-\Lambda_x(u)\right)}e^{b_n\left(\Lambda_{n,x}(u)-\Lambda_x(u)\right)}.
\end{eqnarray*}
Therefore, for every $u>0$,
\begin{eqnarray}
\label{fun}
\sup_{x\in U}\mathbb{P}\left[S_n(x)\geq \delta\right] \leq e^{-b_n\left(u\delta-\sup_{x \in U}\Lambda_x(u)\right)}e^{b_n\sup_{x \in U}\left|\Lambda_{n,x}(u)-\Lambda_x(u)\right|}.
\end{eqnarray}
Similarly, we prove that, for every $u<0$,
\begin{eqnarray}
\label{nrj}
\sup_{x\in U}\mathbb{P}\left[S_n(x)\leq -\delta\right] \leq e^{-b_n\left(-u\delta-\sup_{x \in U}\Lambda_x(u)\right)}e^{b_n\sup_{x \in U}\left|\Lambda_{n,x}(u)-\Lambda_x(u)\right|}.
\end{eqnarray}
The application of Lemma \ref{sup} to (\ref{fun}) and (\ref{nrj}) yields
\begin{eqnarray*}
\sup_{x\in U}\mathbb{P}\left[S_n(x)\geq \delta\right] & \leq & e^{-b_ng_U(\delta)}e^{b_n\sup_{x \in U}\left|\Lambda_{n,x}(\phi(\delta))-\Lambda_x(\phi(\delta))\right|} \\
\sup_{x\in U}\mathbb{P}\left[S_n(x)\leq -\delta\right] & \leq & e^{-b_ng_U(-\delta)}e^{b_n\sup_{x \in U}\left|\Lambda_{n,x}(\phi(-\delta))-\Lambda_x(\phi(-\delta))\right|}
\end{eqnarray*}
and the second part of Lemma \ref{con} provides
\begin{eqnarray*}
\lim_{n \to \infty}\sup_{x \in U}\left|\Lambda_{n,x}(\phi(\delta))-\Lambda_x(\phi(\delta))\right| & = & 0\\
\lim_{n \to \infty}\sup_{x \in U}\left|\Lambda_{n,x}(\phi(-\delta))-\Lambda_x(\phi(-\delta))\right|& = & 0.
\end{eqnarray*}
Consequently, it follows that
\begin{eqnarray*}
\limsup_{n \to \infty}\frac{1}{b_n}\log\sup_{x \in U}\mathbb{P}\left[S_n(x)\geq \delta\right] & \leq & -g_U(\delta)\\
\limsup_{n \to \infty}\frac{1}{b_n}\log\sup_{x \in U}\mathbb{P}\left[S_n(x)\leq -\delta\right] & \leq & -g_U(-\delta)
\end{eqnarray*}
and thus, setting $\tilde{g}_U(\delta)=\min\{g_U(\delta),g_U(-\delta)\}$,
\begin{eqnarray*}
\limsup_{n \to \infty}\frac{1}{b_n}\log\sup_{x \in U}\mathbb{P}\left[\left|S_n(x)\right|\geq \delta\right] & \leq & -\tilde{g}_U(\delta).
\end{eqnarray*}
In order to conclude the proof of Lemma \ref{sort}, let us note that there exists $x_0 \in \overline{U}$ such that $f(x_0)=\|f\|_{U,\infty}$. The application of Proposition \ref{vardp} at the point $x_0$ thus yields
\begin{eqnarray*}
\lim_{n \to \infty}\frac{1}{b_n}\log\mathbb{P}\left[S_n(x_0) \geq \delta\right] & = & -g_U(\delta)\\
\lim_{n \to \infty}\frac{1}{b_n}\log\mathbb{P}\left[S_n(x_0) \leq -\delta\right] & = & -g_U(-\delta)\\
\lim_{n \to \infty}\frac{1}{b_n}\log\mathbb{P}\left[\left|S_n(x_0)\right|\geq \delta\right] & = & -\tilde{g}_U(\delta).
\end{eqnarray*}
The latter relation being due to the straightforward bounds
\begin{eqnarray*}
\max\{\mathbb{P}\left[S_n(x_0) \geq \delta\right] ,\mathbb{P}\left[S_n(x_0) \leq -\delta\right]\} & \leq &  \mathbb{P}\left[\left|S_n(x_0)\right| \geq \delta\right]\\ \mbox{} & & \leq 2\max\{\mathbb{P}\left[S_n(x_0) \geq \delta\right] ,\mathbb{P}\left[S_n(x_0) \leq -\delta\right]\}. \quad\square 
\end{eqnarray*}

\begin{lemma}\label{marem}
Let Assumptions (H1), (H3), (H4)i), (H9)i) and (H10) hold and assume that either $(v_n)\equiv 1$ or (H6) holds.
\begin{enumerate}
\item If $U$ is a bounded set, then, for any $\delta>0$, we have
\begin{eqnarray*}
\lim_{n \to \infty}\frac{v_n^2}{\sum_{i=1}^nh_i^{d+2|\alpha|}}\log\mathbb{P}\left[\sup_{x \in U}v_n\left|\Psi_n^{[\alpha]}(x)\right|\geq \delta\right] & \leq & -\tilde{g}_U(\delta).
\end{eqnarray*}
\item If $U$ is an unbounded set, then, for any $b>0$ and $\delta>0$,
\begin{eqnarray*}
\limsup_{n \to \infty}\frac{v_n^2}{\sum_{i=1}^nh_i^{d+2|\alpha|}}\log\mathbb{P}\left[\sup_{x \in U,\|x\| \leq w_n}v_n\left|\Psi_n^{[\alpha]}(x)\right|\geq \delta\right] & \leq & db-\tilde{g}_U(\delta) 
\end{eqnarray*}
where $w_n=\exp\left(b\sum_{i=1}^nh_i^{d+2|\alpha|}/v_n^2\right)$.
\end{enumerate}

\end{lemma}

\subsubsection*{Proof of Lemma \ref{marem}}
Set $\rho \in ]0,\delta[$, let $\beta$ denote the H\"older order of $\partial^{[\alpha]}K$, and $\|\partial^{[\alpha]}K\|_H$ its corresponding H\"older norm. Set $w_n=\exp\left(b\sum_{i=1}^nh_i^{d+2|\alpha|}/v_n^2\right)$ and 
\begin{eqnarray*}
R_n & = & \left(\frac{\rho n}{2\|\partial^{[\alpha]}K\|_Hv_n\sum_{j=1}^nh_j^{-(d+\beta+|\alpha|)}}\right)^{\frac{1}{\beta}}.
\end{eqnarray*}
We begin with the proof of the second part of  Lemma \ref{marem}. There exist $N'(n)$ points of $\mathbb{R}^d$, $y_1^{(n)},y_2^{(n)},\ldots,y_{N'(n)}^{(n)}$ such that the ball $\{x \in \mathbb{R}^d; \|x\|\leq w_n\}$ can be covered by the $N'(n)$ balls $B_i^{(n)}=\{x\in \mathbb{R}^d; \|x-y_i^{(n)}\|\leq R_n\}$ and such that $\displaystyle N'(n)\leq 2\left(\frac{2w_n}{R_n}\right)$. Considering only the $N(n)$ balls that intersect $\{x \in U;\|x\|\leq w_n\}$, we can write
\begin{eqnarray*}
\{x \in U; \|x\|\leq w_n\}\subset\cup_{i=1}^{N(n)}B_i^{(n)}.
\end{eqnarray*}
For each $i \in \left\{1,\dots,N(n)\right\}$, set $x_i^{(n)}\in B_i^{(n)}\cap U$. We then have: 
\begin{eqnarray*}
\mathbb{P}\left[ {\textstyle\sup_{x \in U, \|x\|\leq w_n}}v_n\left|\Psi_n^{[\alpha]}(x)\right|\geq\delta  \right] & \leq & 
\sum_{i=1}^{N(n)}\mathbb{P}\left[{\textstyle\sup_{x \in B_i^{(n)}}}v_n\left|\Psi_n^{[\alpha]}(x)\right|\geq  \delta \right] \\ & \leq &
 N(n)\max_{1\leq i \leq N(n)} \mathbb{P}\left[\textstyle{\sup_{x \in B_i^{(n)}}}v_n\left|\Psi_n^{[\alpha]}(x)\right|\geq  \delta \right].  
\end{eqnarray*} 
Now, for any $i \in \left\{1,\ldots,N(n)\right\}$ and any $x \in B_i^{(n)}$,
\begin{eqnarray*}
v_n\left|\Psi_n^{[\alpha]}(x)\right| & \leq &
 v_n\left|\Psi_n^{[\alpha]}(x_i^{(n)})\right| + \frac{v_n}{n}\sum_{j=1}^n\frac{1}{h_j^{d+|\alpha|}}\left|\partial^{[\alpha]}K\left(\frac{x-X_j}{h_j}\right)-\partial^{[\alpha]}K\left(\frac{x_i^{(n)}-X_j}{h_j}\right)\right| \\ \mbox{} & & + 
\frac{v_n}{n}\sum_{j=1}^n\frac{1}{h_j^{d+|\alpha|}}\mathbb{E}\left|\partial^{[\alpha]}K\left(\frac{x-X_j}{h_j}\right)-\partial^{[\alpha]}K\left(\frac{x_i^{(n)}-X_j}{h_j}\right)\right| \\ 
& \leq & 
 v_n\left|\Psi_n^{[\alpha]}(x_i^{(n)})\right| +\frac{2 v_n}{n}\|\partial^{[\alpha]}K\|_H\sum_{j=1}^n\frac{1}{h_j^{d+|\alpha|}}\left(\frac{\|x-x_i^{(n)}\|}{h_j}\right)^{\beta} \\
 & \leq &  
 v_n\left|\Psi_n^{[\alpha]}(x_i^{(n)})\right|+\frac{2v_n}{n}\|\partial^{[\alpha]}K\|_H\sum_{j=1}^nh_j^{-(d+\beta+|\alpha|)}R_n^{\beta} \\ & \leq & 
 v_n\left|\Psi_n^{[\alpha]}(x_i^{(n)})\right| + \rho.    
\end{eqnarray*} 
Hence, we deduce that 
\begin{eqnarray*}
\mathbb{P}\left[ {\textstyle\sup_{x \in U,\|x\|\leq w_n}} v_n\left|\Psi_n^{[\alpha]}(x)\right|\geq\delta\right] & \leq &
 N(n)\max_{1\leq i\leq N(n)}\mathbb{P}\left[ v_n\left|\Psi_n^{[\alpha]}(x_i^{(n)})\right| \geq \delta -\rho \right] \\ & \leq &
 N(n)\sup_{x \in U}\mathbb{P}\left[ v_n\left|\Psi_n^{[\alpha]}(x)\right|\geq\delta -\rho \right].
\end{eqnarray*}
Let us at first assume that 
\begin{eqnarray}
\label{bor}
\limsup_{n \to \infty}\frac{v_n^2}{\sum_{i=1}^nh_i^{d+2|\alpha|}}\log N(n)\leq db.
\end{eqnarray}
The application of Lemma \ref{sort} then yields
\begin{eqnarray*}
\lefteqn{\limsup_{n \to \infty}\frac{v_n^2}{\sum_{i=1}^nh_i^{d+2|\alpha|}}\log \mathbb{P}\left[ {\textstyle\sup_{x \in U,\|x\|\leq w_n}} v_n\left|\Psi_n^{[\alpha]}(x)\right|\geq\delta\right]}\\ & \leq & 
\limsup_{n \to \infty}\frac{v_n^2}{\sum_{i=1}^nh_i^{d+2|\alpha|}}\log N(n)-\tilde{g}_U(\delta-\rho)\\ & \leq & 
db-\tilde{g}_U(\delta-\rho).
\end{eqnarray*}
Since this inequality holds for any $\rho \in]0,\delta[$, Part 2 of Lemma \ref{marem} thus follows from the continuity of $\tilde{g}_U$.\\
Let us now establish Relation (\ref{bor}). By definition of $N(n)$ and $w_n$, we have $\log N(n)\leq \log N'(n)\leq db\sum_{i=1}^nh_i^{d+2|\alpha|}/v_n^2+(d+1)\log 2-d\log R_n$, with 
\begin{eqnarray*}
\lefteqn{\frac{v_n^2}{\sum_{i=1}^nh_i^{d+2|\alpha|}}\log R_n}\\  & = & 
\frac{v_n^2}{\beta\sum_{i=1}^nh_i^{d+2|\alpha|}}\left[\log \rho +\log n -\log\left(2\|\partial^{[\alpha]}K\|_H\right)-\log v_n - \log\left(\sum_{j=1}^nh_j^{-(d+\beta+|\alpha|)}\right)\right],
\end{eqnarray*}
which goes to zero in view of (H10) and \eqref{rel}. Thus, \eqref{bor} is proved, and the proof of part 2 of Lemma \ref{marem} is completed.\\

Let us now consider part 1 of Lemma \ref{marem}. This part is proved by following the same steps as for part 2, except that the number $N(n)$ of balls covering $U$ is at most the integer part of $\left(\Delta/R_n\right)^d$, where $\Delta$ denotes the diameter of $\overline{U}$. Relation (\ref{bor}) then becomes
\begin{eqnarray*}
\limsup_{n \to \infty}\frac{v_n^2}{\sum_{i=1}^nh_i^{d+2|\alpha|}}\log R_n\leq 0
\end{eqnarray*}
and Lemma \ref{marem} is proved. $\square$
\begin{lemma}\label{moussa}
Let Assumptions (H1), (H3), (H4)i) and (H11)ii) hold. Assume that either $(v_n)\equiv 1$ or (H6), (H10) and (H11)i) hold.
Moreover assume that $\partial^{[\alpha]}f$ is continuous.
For any $b>0$ if we set $w_n=\exp\left(b\sum_{i=1}^nh_i^{d+2|\alpha|}/v_n^2\right)$ then, for any $\rho>0$, we have, for $n$ large enough,
\begin{eqnarray*}
\sup_{x\in U,\|x\|\geq w_n}\frac{v_n}{n}\sum_{i=1}^n\frac{1}{h_i^{d+|\alpha|}}\left|\mathbb{E}\left[\partial^{[\alpha]}K\left(\frac{x-X_i}{h_i}\right)\right]\right|\leq \rho.
\end{eqnarray*}
\end{lemma}

\subsubsection*{Proof of Lemma \ref{moussa}}
We have
\begin{eqnarray}
\label{pour}
\frac{v_n}{n}\sum_{i=1}^n\frac{1}{h_i^{d+|\alpha|}}\mathbb{E}\left[\partial^{[\alpha]}K\left(\frac{x-X_i}{h_i}\right)\right] & = & \frac{v_n}{n}\sum_{i=1}^n\int_{\mathbb{R}^d}K(z)\partial^{[\alpha]}f(x-h_iz)dz.
\end{eqnarray}
Set $\rho>0$. In the case $(v_n)\equiv 1$, set $M$ such that $\|\partial^{[\alpha]}f\|_{\infty}\int_{\|z\|> M}\left|K(z)\right|dz \leq\rho/2$ ; we have
\begin{eqnarray*}
\lefteqn{\frac{v_n}{n}\sum_{i=1}^n\frac{1}{h_i^{d+|\alpha|}}\left|\mathbb{E}\left[\partial^{[\alpha]}K\left(\frac{x-X_i}{h_i}\right)\right]\right|} \\  & \leq & 
\frac{\rho}{2} +  \partial^{[\alpha]}f(x)\int_{\|z\|\leq M}\left|K(z)\right|dz +\frac{1}{n}\sum_{i=1}^n\int_{\|z\|\leq M}\left|K(z)\right| \left| \partial^{[\alpha]}f(x-h_iz)-\partial^{[\alpha]}f(x)\right|dz.
\end{eqnarray*}
Lemma \ref{moussa} then follows from the fact that $\partial^{[\alpha]}f$ fulfills (H11)ii). As a matter of fact, this condition implies that $\lim_{\|x\| \to \infty,x \in\overline{U}}\partial^{[\alpha]}f(x)=0$ and that the third term in the right-hand-side of the previous inequality goes to $0$ as $n \to \infty$ (by the dominated convergence).\\
Let us now assume that $\lim_{n \to \infty}v_n=\infty$; relation (\ref{pour}) can be rewritten as
\begin{eqnarray*}
\frac{v_n}{n}\sum_{i=1}^n\frac{1}{h_i^{d+|\alpha|}}\mathbb{E}\left[\partial^{[\alpha]}K\left(\frac{x-X_i}{h_i}\right)\right]  & = & 
\frac{v_n}{n}\sum_{i=1}^n\int_{\|z\|\leq w_n/2}K(z)\partial^{[\alpha]}f(x-h_iz)dz \\ \mbox{} & & +\frac{v_n}{n}\sum_{i=1}^n\int_{\|z\|> w_n/2}K(z)\partial^{[\alpha]}f(x-h_iz)dz.
\end{eqnarray*}
Set $\rho>0$; on the one hand, we have
\begin{eqnarray*}
\|x\| \geq w_n \quad \mbox{and}\quad \|z\|\leq w_n/2 & \Rightarrow & \|x-h_iz\|\geq w_n\left(1-h_i/2\right)\\
& \Rightarrow & \|x-h_iz\| \geq w_n/2 \quad \mbox{for $n$ large enough}.
\end{eqnarray*}
Set $M_f=\sup_{x \in \mathbb{R}^d}\|x\|^{\eta}\partial^{[\alpha]}f(x)$. Assumption (H11)ii) implies that, for $n$ sufficiently large, 
\begin{eqnarray*}
\sup_{\|x\|\geq w_n}\frac{v_n}{n}\sum_{i=1}^n\int_{\|z\|\leq w_n/2}\left|K(z)\partial^{[\alpha]}f(x-h_iz)\right|dz &\leq & 
\sup_{\|x\|\geq w_n}\frac{v_n}{n}\sum_{i=1}^n\int_{\|z\|\leq w_n/2}\left|K(z)\right|M_f\|x-h_iz\|^{-\eta}dz \\ & \leq & 
2^{\eta}\frac{v_n}{w_n^{\eta}}M_f\int_{\mathbb{R}^d}\left|K(z)\right|dz \\ & \leq & 
\frac{\rho}{2}.
\end{eqnarray*}
On the other hand, we note that, in view of assumptions (H10) and (H11)i),
\begin{eqnarray*}
\sup_{\|x\|\geq w_n}\frac{v_n}{n}\sum_{i=1}^n\int_{\|z\|> w_n/2}\left|K(z)\partial^{[\alpha]}f(x-h_iz)\right|dz  \leq 
2^{\zeta}\frac{v_n}{w_n^{\zeta}}\|\partial^{[\alpha]}f\|_{\infty}\int_{\|z\|> w_n/2}\|z\|^{\zeta}\left|K(z)\right|dz \leq \frac{\rho}{2}
\end{eqnarray*}
(for $n$ large enough). As a matter of fact, we have by assumptions (H6) and (H10), $\forall \xi >0$
\begin{eqnarray*}
\frac{v_n}{w_n^{\xi}} = \exp\left\{-\xi\log w_n\left(1-\frac{\log v_n}{\xi \log w_n}\right)\right\}\stackrel{n \to \infty}\longrightarrow0.
\end{eqnarray*}
This concludes the proof of Lemma \ref{moussa}. $\square$\\

Since $\partial^{[\alpha]}K$ is a bounded function that vanishes at infinity, we have $\lim_{\|x\| \to \infty}|\Psi_n^{[\alpha]}(x)|=0$ for every given $n\geq 1$. Moreover, since $\partial^{[\alpha]}K$  is assumed to be continuous, $\Psi_n^{[\alpha]}$ is continuous, and this ensures the existence of a random variable $s_n$ such that 
\begin{eqnarray*}
\left|\Psi_n^{[\alpha]}(s_n)\right| & = & \sup_{x \in U}\left|\Psi_n^{[\alpha]}(x)\right|.
\end{eqnarray*}

\begin{lemma}\label{rougui}
Let Assumptions (H1), (H3), (H4)i), (H8)i), (H9)ii) and (H10) hold. Suppose either $(v_n)\equiv 1$ or (H6) and (H11) hold.
For any $b>0$, set $w_n=\exp\left(b\sum_{i=1}^nh_i^{d+2|\alpha|}/v_n^2\right)$; then, for any $\delta >0$, we have
\begin{eqnarray*}
\limsup_{n \to \infty}\frac{v_n^2}{\sum_{i=1}^nh_i^{d+2|\alpha|}}\log\mathbb{P}\left[\|s_n\| \geq w_n \quad \mbox{and}\quad v_n\left|\Psi_n^{[\alpha]}(s_n)\right|\geq \delta\right] & \leq & -b\xi.
\end{eqnarray*}

\end{lemma}

\subsubsection*{Proof of Lemma \ref{rougui}}
We first note that $s_n \in \overline{U}$ and therefore
\begin{eqnarray*}
\lefteqn{\|s_n\|\geq w_n \ \ \mbox{and}\ \ 
v_n\left|\Psi_n^{[\alpha]}(s_n)\right| \geq \delta}\\
& \Rightarrow & \|s_n\|\geq w_n \ \ \mbox{and}\ \
\frac{v_n}{n}\left|\sum_{i=1}^n\frac{1}{h_i^{d+|\alpha|}}\partial^{[\alpha]}K\left(\frac{s_n-X_i}{h_i}\right)\right|+\frac{v_n}{n}\mathbb{E}\left|\sum_{i=1}^n\frac{1}{h_i^{d+|\alpha|}}\partial^{[\alpha]}K\left(\frac{s_n-X_i}{h_i}\right)\right| \geq \delta\\
& \Rightarrow & \|s_n\|\geq w_n \ \ \mbox{and}\ \ 
\frac{v_n}{n}\sum_{i=1}^n\frac{1}{h_i^{d+|\alpha|}}\left|\partial^{[\alpha]}K\left(\frac{s_n-X_i}{h_i}\right)\right|\geq \delta \mbox{}\\ & & -   \sup_{\|x\| \geq w_n, x\in \overline{U}}\frac{v_n}{n}\sum_{i=1}^n\frac{1}{h_i^{d+|\alpha|}}\mathbb{E}\left|\partial^{[\alpha]}K\left(\frac{x-X_i}{h_i}\right)\right|\\
\end{eqnarray*}
Set $\rho \in ]0,\delta[$; the application of Lemma \ref{moussa} ensures that, for $n$ large enough,
\begin{eqnarray*}
\lefteqn{\|s_n\|\geq w_n \ \ \mbox{and}\ \ 
v_n\left|\Psi_n^{[\alpha]}(s_n)\right| \geq \delta}\\
& \Rightarrow & \|s_n\|\geq w_n \ \ \mbox{and}\ \
\frac{v_n}{n}\sum_{i=1}^n\frac{1}{h_i^{d+|\alpha|}}\left|\partial^{[\alpha]}K\left(\frac{s_n-X_i}{h_i}\right)\right|\geq \delta -\rho.
\end{eqnarray*}
Set $\displaystyle\kappa =\sup_{x \in \mathbb{R}^d}\|x\|^{\gamma}|\partial^{[\alpha]}K(x)|$ (see Assumption (H9)ii)). We obtain, for $n$ sufficiently large, 
\begin{eqnarray*}
\lefteqn{\|s_n\|\geq w_n \ \ \mbox{and}\ \ 
v_n\left|\Psi_n^{[\alpha]}(s_n)\right|
\geq \delta}\\
& \Rightarrow & \|s_n\|\geq w_n \ \ \mbox{and}\ \ 
\exists i \in \left\{1, \ldots, n\right\}\ \ \mbox{such that}\ \  
\frac{v_n}{h_i^{d+|\alpha|}}\left|\partial^{[\alpha]}
K\left(\frac{s_n-X_i}{h_i}\right)\right| \geq \delta - \rho\\  
& \Rightarrow & \|s_n\|\geq w_n \ \ \mbox{and}\ \ 
\exists i \in \left\{1, \ldots, n\right\} \ \ \mbox{such that}\ \  
\kappa h_i^{\gamma}\geq \frac{h_i^{d+|\alpha|}}{v_n}\|s_n-X_i\|^{\gamma}
(\delta-\rho) \\
& \Rightarrow & \|s_n\|\geq w_n \ \ \mbox{and}\ \ 
\exists i \in \left\{1, \ldots, n\right\} \ \ \mbox{such that}\ \ 
\big|\|s_n\|-\|X_i\|\big|\leq 
\Bigg[\frac{\kappa v_nh_i^{\gamma-d-|\alpha|}}{\delta-\rho}
\Bigg]^{\frac{1}{\gamma}}\\
& \Rightarrow & \|s_n\| \geq w_n \ \ \mbox{and}\ \ 
\exists i \in \left\{1, \ldots, n\right\} \ \ \mbox{such that}\ \  
\|X_i\|\geq \|s_n\|-\Bigg[\frac{\kappa v_nh_i^{\gamma-d-|\alpha|}}
{\delta-\rho}\Bigg]^{\frac{1}{\gamma}}\\
& \Rightarrow & \exists i \in \left\{1, \ldots, n\right\}
\ \ \mbox{such that}\ \ \|X_i\|\geq w_n\left(1-u_{n,i}\right)\quad \mbox{with} \quad u_{n,i} = w_n^{-1}v_n^{\frac{1}{\gamma}}h_i^{\frac{\gamma-d-|\alpha|}{\gamma}}\left(\frac{\kappa}{\delta-\rho}\right)^{\frac{1}{\gamma}}.
\end{eqnarray*}
Assume for the moment that 
\begin{eqnarray}
\label{placebo}
\lim_{n \to \infty}u_{n,i}=0.
\end{eqnarray}
It then follows that $1-u_{n,i}>0$ for $n$ sufficiently large; therefore we can deduce that (see Assumption (H8)i)):
\begin{eqnarray*}
\mathbb{P}\left[\|s_n\| \geq w_n \quad \mbox{and}\quad v_n\left|\Psi_n^{[\alpha]}(s_n)\right|\geq \delta\right] &\leq & 
\sum_{i=1}^n\mathbb{P}\left[\|X_i\|^{\xi}\geq w_n^{\xi}\left(1-u_{n,i}\right)^{\xi}\right]\\ & \leq &
 \sum_{i=1}^n\mathbb{E}\left(\|X_i\|^{\xi}\right)w_n^{-\xi}\left(1-u_{n,i}\right)^{-\xi}\\ & \leq & 
n\mathbb{E}\left(\|X_1\|^{\xi}\right)w_n^{-\xi}\max_{1\leq i \leq n}\left(1-u_{n,i}\right)^{-\xi}.
\end{eqnarray*}
Consequently,
\begin{eqnarray*}
\lefteqn{\frac{v_n^2}{\sum_{i=1}^nh_i^{d+2|\alpha|}}\log \mathbb{P}\left[\|s_n\| \geq w_n \quad \mbox{and}\quad v_n\left|\Psi_n^{[\alpha]}(s_n)\right|\geq \delta\right]}\\  &\leq & 
\frac{v_n^2}{\sum_{i=1}^nh_i^{d+2|\alpha|}}\left[\log n+\log \mathbb{E}\left(\|X_1\|^{\xi}\right)-b\frac{\sum_{i=1}^nh_i^{d+2|\alpha|}}{v_n^2}-\xi\log\max_{1\leq i\leq n}\left(1-u_{n,i}\right)\right],
\end{eqnarray*}
and, thanks to assumption (H10), it follows that
\begin{eqnarray*}
\limsup_{n \to \infty}\frac{v_n^2}{\sum_{i=1}^nh_i^{d+2|\alpha|}}\log \mathbb{P}\left[\|s_n\| \geq w_n \quad \mbox{and}\quad v_n\left|\Psi_n^{[\alpha]}(s_n)\right|\geq \delta\right] & \leq & -b\xi. 
\end{eqnarray*}
Let us now prove relation (\ref{placebo}). We expand
\begin{eqnarray*}
u_{n,i} & = & \exp\left(-b\frac{\sum_{i=1}^nh_i^{d+2|\alpha|}}{v_n^2}\left[1-\frac{1}{b\gamma}\frac{v_n^2\log v_n}{\sum_{i=1}^nh_i^{d+2|\alpha|}}-\frac{\gamma-d-|\alpha|}{b\gamma}\frac{v_n^2\log\left(h_i\right)}{\sum_{i=1}^nh_i^{d+2|\alpha|}}\right]\right)\left(\frac{\kappa}{\delta-\rho}\right)^{\frac{1}{\gamma}}
\end{eqnarray*}
and assumptions (H6) and (H10) ensure that $\lim_{n \to \infty}u_{n,i}=0$ and thus Lemma \ref{rougui} is proved. $\square$\\
\subsubsection*{Proof of Proposition \ref{univardev}}
Let us at first note that the lower bound
\begin{eqnarray}
\label{sone}
\liminf_{n \to \infty}\frac{v_n^2}{\sum_{i=1}^nh_i^{d+2|\alpha|}}\log\mathbb{P}\left[\sup_{x \in U}v_n\left|\Psi_n^{[\alpha]}(x)\right|\geq \delta\right] \geq -\tilde{g}_U(\delta)
\end{eqnarray}
follows from the application of Proposition \ref{vardp} at a point $x_0\in\overline{U}$ such that $f(x_0)=\|f\|_{U,\infty}$.\\
In the case $U$ is bounded, Proposition \ref{univardev} is thus a straightforward consequence of (\ref{sone}) and of the first part of Lemma \ref{marem}. Let us now consider the case $U$ is unbounded.\\
Set $\delta>0$ and, for any $b>0$ set $w_n=\exp\left(b\sum_{i=1}^nh_i^{d+2|\alpha|}/v_n^2\right)$. Since, by definition of $s_n$,
\begin{eqnarray*}
\lefteqn{\mathbb{P}\left[\sup_{x \in U}v_n\left|\Psi_n^{[\alpha]}(x)\right|\geq \delta\right]}\\ & \leq & 
\mathbb{P}\left[\sup_{x \in U,\|x\| \leq w_n}v_n\left|\Psi_n^{[\alpha]}(x)\right|\geq \delta\right] + \mathbb{P}\left[\|s_n\| \geq w_n \quad \mbox{and} \quad v_n\left|\Psi_n^{[\alpha]}(s_n)\right|\geq \delta\right]
\end{eqnarray*}
it follows from Lemmas \ref{marem} and \ref{rougui} that
\begin{eqnarray*}
\limsup_{n \to \infty}\frac{v_n^2}{\sum_{i=1}^nh_i^{d+2|\alpha|}}\log\mathbb{P}\left[\sup_{x \in U}v_n\left|\Psi_n^{[\alpha]}(x)\right|\geq \delta\right] & \leq & \max\{-b\xi;db-\tilde{g}_U(\delta)\}
\end{eqnarray*}
and consequently
\begin{eqnarray*}
\limsup_{n \to \infty}\frac{v_n^2}{\sum_{i=1}^nh_i^{d+2|\alpha|}}\log\mathbb{P}\left[\sup_{x \in U}v_n\left|\Psi_n^{[\alpha]}(x)\right|\geq \delta\right] & \leq & \inf_{b>0}\max\{-b\xi;db-\tilde{g}_U(\delta)\}.
\end{eqnarray*}
Since the infimum in the right-hand-side of the previous bound is achieved for $b=\tilde{g}_U(\delta)/\left(\xi+d\right)$ and equals $-\xi\tilde{g}_U(\delta)/\left(\xi+d\right)$, we obtain the upper bound
\begin{eqnarray*}
\limsup_{n \to \infty}\frac{v_n^2}{\sum_{i=1}^nh_i^{d+2|\alpha|}}\log\mathbb{P}\left[\sup_{x \in U}v_n\left|\Psi_n^{[\alpha]}(x)\right|\geq \delta\right] & \leq& -\frac{\xi}{\xi+d}\tilde{g}_U(\delta)
\end{eqnarray*}
which concludes the proof of Proposition \ref{univardev}. $\square$
\subsection{Proof of Proposition \ref{bia}}
Let us set $g=\partial^{[\alpha]}f$, $D^jg$ $\big(j \in \{1,\ldots,q\}\big)$ the $j$-th differential of $g$,
$y = (y_1, \dots ,y_d)\in \mathbb{R}^d$ and $y^{(j)}=(y, \ldots, y) \in (\mathbb{R}^d)^j$. With these notations, 
$$D^jg(x)(y^{(j)})=\sum_{\alpha_1+\cdots+\alpha_d=j}\frac{\partial^jg}{\partial y_1^{\alpha_1}\ldots \partial y_d^{\alpha_d}}(x)y_1^{\alpha_1}\ldots y_d^{\alpha_d}.$$
By successive integrations by parts (and using the fact that the partial derivatives of $K$ vanish at infinity, see Assumptions (H4)i)), we have 
\begin{eqnarray*}
\mathbb{E}\Big[\partial^{[\alpha]}f_n(x)\Big] & = & \frac{1}{n}\sum_{i=1}^n\frac{1}{h_i^{d+|\alpha|}}\mathbb{E}\Bigg[\partial^{[\alpha]}K\left(\frac{x-X_i}{h_i}\right)\Bigg] \\ & = & \frac{1}{n}\sum_{i=1}^n\frac{1}{h_i^{d+|\alpha|}}\int_{\mathbb{R}^d}\partial^{[\alpha]}K\left(\frac{x-y}{h_i}\right)f(y)dy \\ & = & \frac{1}{n}\sum_{i=1}^n\frac{1}{h_i^d}\int_{\mathbb{R}^d}K\left(\frac{x-y}{h_i}\right)g(y)dy \\ & = & \frac{1}{n}\sum_{i=1}^n \int_{\mathbb{R}^d}K(y)g(x-h_iy)dy.
\end{eqnarray*}
Hence, using assumption (H7)i) and the fact that $\partial^{[\alpha]}f$ is $q$-times differentiable, it comes 
\begin{eqnarray}
\label{eqd}
\lefteqn{\mathbb{E}\Big[\partial^{[\alpha]}f_n(x)\Big]-\partial^{[\alpha]}f(x)}  \nonumber\\   & = &  
\frac{1}{n}\sum_{i=1}^n\int_{\mathbb{R}^d}K(y)\big[g(x-h_iy)-g(x)\big]dy \nonumber \\ & = & 
\frac{1}{n}\sum_{i=1}^nh_i^q\int_{\mathbb{R}^d}K(y)\Bigg[\frac{g(x-h_iy)-g(x)-\sum_{j=1}^{q-1}\frac{(-1)^j}{j!}h_i^jD^jg(x)(y^{(j)})}{h_i^q}\Bigg]dy.
\end{eqnarray}
Let us set 
\begin{eqnarray*}
U_i(x)  & = & \int_{\mathbb{R}^d}K(y)\left[\frac{g(x-h_iy)-g(x)-\sum_{j=1}^{q-1}\frac{(-1)^j}{j!}h_i^jD^jg(x)(y^{(j)})}{h_i^q}\right]dy\quad  \mbox{and}\\ U_{\infty}(x) & = &  \frac{(-1)^q}{q!}\int_{\mathbb{R}^d}D^qg(x)(y^{(q)})K(y)dy.  
\end{eqnarray*}
 We clearly have 
\begin{eqnarray}
\label{lim}
\lim_{i \to \infty}U_i(x) & = & U_{\infty}(x)
\end{eqnarray}
and therefore, $\forall \varepsilon >0$, $\exists i_0 \in \mathbb{R}$ such that $\forall i\geq i_0$, $\left|U_i(x)-U_{\infty}(x)\right|\leq \varepsilon$. 
\begin{itemize}
\item If $\sum_i h_i^q=\infty$, then 
\begin{eqnarray*}
\lefteqn{\left|\frac{n}{\sum_{i=1}^nh_i^q}\left[\mathbb{E}\left(\partial^{[\alpha]}f_n(x)\right)-\partial^{[\alpha]}f(x)\right] -U_{\infty}(x)\right|}\\ & = &
\left|\frac{\sum_{i=1}^nh_i^qU_i(x)}{\sum_{i=1}^nh_i^q}-U_{\infty}(x) \right| \\ & \leq &
\frac{\sum_{i=1}^{i_0-1}h_i^q\left|U_i(x)-U_{\infty}(x)\right|+\sum_{i=i_0}^nh_i^q\left|U_i(x)-U_{\infty}(x)\right|}{\sum_{i=1}^nh_i^q} \\ & \leq & 
2\varepsilon. 
\end{eqnarray*}
\item If $\sum_ih_i^q<\infty$, we can write
\begin{eqnarray*}
\frac{n}{\sum_{i=1}^nh_i^q}\left[\mathbb{E}\left(\partial^{[\alpha]}f_n(x)\right)-\partial^{[\alpha]}f(x)\right] & = &
\frac{\sum_{i=1}^nh_i^qU_i(x)}{\sum_{i=1}^nh_i^q}. 
\end{eqnarray*}
\end{itemize}
In view of (\ref{lim}), for $x$ fixed and for all $i \in \mathbb{N}$, the sequence $(U_i(x))_i$ is bounded and thus Part 1 of Proposition \ref{bia} is completed. Let us now prove Part 2. \\
Since the bracketed term in (\ref{eqd}) is bounded by $\sup_{x\in\mathbb{R}^d}\|D^qg(x)\|=M_q$ (see Assumption (H7)iii)), Part 2 follows. $\square$ \\

\subsection{Proof of Proposition \ref{rate}}
\begin{itemize}
 \item Since $\left|e^t-1\right| \leq \left|t\right|e^{\left|t\right|}$ $\forall t \in \mathbb{R}$, and thanks to the boundedness and integrability of $K$, we have
\begin{eqnarray*}
\int_{[0,1]\times\mathbb{R}^d}s^{-ad}\left|e^{s^{ad}\frac{u}{1-ad}K(z)}-1\right|dsdz & \leq & \frac{|u|}{1-ad}e^{\frac{|u|}{1-ad}\|K\|_{\infty}}\int_{[0,1]\times\mathbb{R}^d}s^{-ad}\left|K(z)\right|dsdz<\infty
\end{eqnarray*}
which ensures the existence of $\psi$.
It is straightforward to check that $\psi$ is twice differentiable, with
\begin{eqnarray*}
\psi'(u) & = & \frac{1}{1-ad}\int_{[0,1]\times\mathbb{R}^d}K(z)e^{s^{ad}\frac{u}{1-ad}K(z)}dsdz, \\  \psi''(u) & = & \frac{1}{\left(1-ad\right)^2}\int_{[0,1]\times\mathbb{R}^d}s^{ad}\left(K(z)\right)^2e^{s^{ad}\frac{u}{1-ad}K(z)}dsdz. 
\end{eqnarray*}
Since $\psi''(u)>0$ $\forall u \in \mathbb{R}$, $\psi'$ is increasing on $\mathbb{R}$, and $\psi$ is strictly convex on $\mathbb{R}$. It follows that its Cramer transform $I$ is a good rate function on $\mathbb{R}$ (see Dembo and Zeitouni (1998)) and $(i)$ of Proposition \ref{rate} is proved.
\item Let us now assume that $\lambda (S_-)=0$. We then have
\begin{eqnarray*}
\lim_{u \to -\infty}\psi'(u)=0 \quad \mbox{and}\quad \lim_{u \to +\infty}\psi'(u)=+\infty,
\end{eqnarray*}
so that the range of $\psi'$ is $]0,+\infty[$. Moreover $\lim_{u \to -\infty}\psi(u)=-\lambda(S_+)/(1-ad)$ (which can be $-\infty$). This implies in particular that $I(0)=\lambda(S_+)/(1-ad)$. Now, when $t<0$, $\lim_{u \to -\infty}\left(ut-\psi(u)\right)=+\infty$, and $I(t)=+\infty$. Since $\psi'$ is increasing with range $]0,+\infty[$, when $t>0$, $\sup_u(ut-\psi(u))$ is reached for $u_0(t)$ such that $\psi'(u_0(t))=t$, i.e. for $u_0(t)=(\psi')^{-1}(t)$; this prove \eqref{mel}. (Note that, since $\psi''(t)>0$, the function $t\mapsto u_0(t)$ is differentiable on $]0,+\infty[$). Now, differentiating \eqref{mel}, we have  
\begin{eqnarray*}
I'(t) & = & u_0(t)+tu_0'(t)-\psi'(u_0(t))u_0'(t) \\ 
& = & (\psi')^{-1}(t)+tu_0'(t)-tu_0'(t) \\
& = & (\psi')^{-1}(t).
\end{eqnarray*}
Since $(\psi')^{-1}$ is an increasing function on $]0,+\infty[$, it follows that $I$ is strictly convex on $]0,+\infty[$ (and differentiable). Thus $(ii)$ is proved.\\
Now, since $\lambda(S_-)=0$, $\psi'(0)=1/(1-ad)>0$; we have 
\begin{eqnarray*}
I'(t)=0 \quad \Leftrightarrow \quad(\psi')^{-1}(t)=0 \quad \Leftrightarrow\quad \psi'(0)=t\quad\Leftrightarrow \quad t=1/(1-ad).
\end{eqnarray*}
Then $I'\left(1/(1-ad)\right)=0$, and $I\left(1/(1-ad)\right)=0$ is the unique global minimum of $I$ on $]0,+\infty[$. This proves $(iv)$ when $\lambda(S_-)=0$.
\item Assume that $\lambda(S_-)>0$. In this case, $\psi'$ can be rewritten as
\begin{eqnarray*}
\psi'(u)
 & = & \frac{1}{1-ad}\int_{[0,1]\times\left(\mathbb{R}^d\cap S_+\right)}K(z)e^{s^{ad}\frac{u}{1-ad}K(z)}dsdz \\ \mbox{} && + \frac{1}{1-ad}\int_{[0,1]\times\left(\mathbb{R}^d\cap S_-\right)}K(z)e^{s^{ad}\frac{u}{1-ad}K(z)}dsdz
\end{eqnarray*}
and we have 
\begin{eqnarray*}
\lim_{u \to -\infty}\psi'(u)=-\infty \quad \mbox{and}\quad \lim_{u \to +\infty}\psi'(u)=+\infty 
\end{eqnarray*}
so that the range of $\psi'$ is $\mathbb{R}$ in this case. The proof of $(iii)$ and the case $\lambda(S_-)>0$ of $(iv)$ follows the same lines as previously, except that, in the present case, $(\psi')^{-1}$ is defined on  $\mathbb{R}$, and not only on $]0,+\infty[$. $\square$
\end{itemize}

\end{document}